\documentclass[11pt]{article}
\usepackage{amssymb}
\newtheorem{cor}{Corollary}[section]
\newtheorem{pro}{Proposition}[section]
\newtheorem{thm}{Theorem}[section]

\makeatletter\@addtoreset{equation}{section}
\renewcommand\theequation{\thesection.\@arabic\c@equation}

\newenvironment{Literature}[1]
{\begin{thebibliography}{#1}}{\end{thebibliography}}

\begin{document}

\begin{center}
{\LARGE \bf  Lower Bounds of the Dirac Eigenvalues }     \\

\bigskip
{\LARGE   \bf  on Compact  Riemannian Spin Manifolds }       \\

\bigskip
{\LARGE  \bf with Locally Product Structure  }

\bigskip  \bigskip
{\large  Eui Chul Kim}
\end{center}

\bigskip   \noindent
Department of Mathematics Education,  Andong National University,  \\
388 Songchon-dong, Andong,  Kyoungsangbuk-do, 760-749, South  Korea     \\
e-mail: eckim@andong.ac.kr

\bigskip \bigskip \noindent
{\bf Abstract.}   We study some similarities between almost
product Riemannian structures and almost Hermitian structures.
Inspired by the similarities, we prove lower eigenvalue estimates
for the Dirac operator on compact  Riemannian spin manifolds with
locally product structures.  We also provide some examples
(limiting manifolds) for the limiting case of the estimates.


\bigskip  \bigskip \bigskip
\noindent
\section{Introduction}

\bigskip  \noindent
Let $(M^n, g)$ be an n-dimensional (smooth oriented) Riemannian
spin manifold without boundary. When $M^n$ is compact, the
spectrum of the Dirac operator $D$ of $(M^n, g)$ is discrete and
real. The first sharp estimate for the smallest absolute value of
eigenvalues $\lambda$ of the Dirac operator $D$ was obtained by
Friedrich [8].  Using a suitable deformation of Riemannian spin
connection, he proved the inequality
\begin{equation}  \lambda^2  \  \geq  \  \frac{n}{4(n-1)} \inf_{M} S     \end{equation}
on manifolds $(M^n, g)$ with positive scalar curvature $S > 0$.
Equality occurs if and only if  $(M^n , g )$ admits a (real) {\it
Killing spinor} , i.e.,  a non-trivial solution $\psi$ to the {\it
Killing  equation} for spinor fields,
\[   \nabla_X   \psi  =  - \, \frac{\lambda}{n} \,  X  \cdot  \psi  ,   \qquad  \lambda  \not= 0  \in {\mathbb R} ,  \]
where $X$ is an arbitrary vector field on $M^n$ and the dot $"
\cdot "$ indicates the Clifford multiplication [3, 10].
 The inequality (1.1) has been improved in several directions [5, 6, 9, 11, 12, 13, 14, 16, 17, 18, 19].
We refer to [7] for  a concise exposition of the first eigenvalue
(the smallest absolute value) estimates as well as the
classification problems of limiting manifolds.

\bigskip
An optimal lower bound of the Dirac eigenvalues depends on the
geometric structure (the holonomy group) that the considered
manifold may possess. A remarkable observation was made by
 Hijazi  that, if $(M^n, g)$ admits a parallel k-form, $0 < k < n$, then there exist no  Killing spinors (see [3], p.32).
Furthermore,  if $(M^n, g)$ possesses a locally product structure,
then there exist no Killing spinors (see [3], p.35). An
interesting improvement of the Friedrich inequality (1.1)  in
these directions was found by Alexandrov, Grantcharov and Ivanov
[1, 20]. They proved that, if $(M^n, g), \, n \geq 3,$ admits a
parallel 1-form, then any eigenvalue $\lambda$ of the Dirac
operator  $D$ satisfies
\begin{equation}  \lambda^2  \  \geq  \  \frac{n-1}{4(n-2)} \inf_{M} S .  \end{equation}
(Equality occurs  if and only if there exists a non-trivial
solution to  the field equation (1.4) below,  with $n_1 = n-1$ and
$n_2 =1$).

\bigskip
In this paper we study  some similarities between {\it almost
product Riemannian structures}  and {\it almost Hermitian
structures}.  In Section 2  we translate some basic results  in
K\"{a}hler spin geometry [13, 14, 15] into the forms to be
appropriate for almost product Riemannian manifolds.  Inspired by
the similarities,   we will prove lower eigenvalue estimates for
the Dirac operator on compact  Riemannian spin manifolds with
locally product structure (Theorem 1.1 and 1.2). Our new
inequalities contain the inequality (1.2) as a special case. To
state the main result of the paper precisely, we now recall some
basic facts from almost product Riemannian geometry [23].

\bigskip  \noindent
A Riemannian manifold $(M^n, g)$ is called  {\it locally decomposable}  if it admits a (1,1)-tensor field $\phi$ with the following properties :  \\
(i)  \,  $\phi^2 (X)  =  X $  and \, $g (  \phi X,   \phi  Y) =  g
(X,   Y) $  \,
for all vector fields $X, Y$.    \\
(ii)  \, $\nabla \phi = 0$, i.e., $\phi$  is parallel with respect to the Levi-Civita connection $\nabla$.   \\
In case that  $(M^n , g,  \phi )$ is  locally decomposable,  the
tangent bundle $T(M^n)$ decomposes into  $T(M^n) = T^+(M^n)
\oplus T^-(M^n)$ under the action of the endomorphism $\phi$,
where
\[   T^{\pm}(M^n) :=  \{ \,  Z \in T(M^n) :  \phi(Z) =  \pm Z   \,  \}  .   \]
Due to the condition $\nabla \phi = 0$, the distributions $
T^{\pm}(M^n)$ are integrable.  In fact, around each point $x \in
M^n$, there is an open neighbourhood $U$  that has a Riemannian
product structure of the form $(U, g)=( U_1 \times U_2,  \, g_1 +
g_2)$. If $(M^n, g)$ is simply connected and complete, then there
is a global splitting $(M^n , g) = (M_1 \times M_2 , \, g_1 + g_2
)$ (see [21], p.228).

\bigskip
Let  $( M^n, g, \phi)$ be a locally decomposable Riemannian
manifold with a fixed spin structure, and let $(E_1, \ldots, E_n)$
be a local orthonormal frame field. Then the spin derivative
$\nabla $ and the Dirac operator $D$ of $(M^n, g, \phi)$, acting
on sections $ \psi  \in \Gamma(\Sigma(M^n))$ of the spinor bundle
$\Sigma(M^n)$, are locally expressed as
\[             \nabla_X   \psi  = X(\psi) + \frac{1}{4} \sum_{i=1}^n
E_i \cdot \nabla_X  E_i \cdot \psi         \] and
\[    D \psi = \sum_{i=1}^n E_i \cdot \nabla_{E_i} \psi,     \]
respectively. As in the K\"{a}hlerian case [13], let us define the
{\it twist } $\widetilde{D}$ of the Dirac operator  $D$ by
\[
\widetilde{D} \psi   =  \sum_{i=1}^n   E_i \cdot \nabla_{\phi
(E_i)} \psi   = \sum_{j=1}^{n} \phi(E_j) \cdot \nabla_{E_j} \psi .
\]

\bigskip  \noindent
In Section 3 we will prove the following theorems.

\noindent
\begin{thm}
Let  $( M^n, g, \phi), \, n \geq 4,$ be a locally decomposable
Riemannian spin manifold with positive scalar curvature $S > 0$.
Assume that $( M^n, g, \phi)$ is compact and the dimension $n_1$
of $T^+(M^n)$ is equal to the dimension $n_2$ of $T^-(M^n)$ (i.e.,
$n= n_1 + n_2 = 2 n_1$). Then any eigenvalue $\lambda$ of the
Dirac operator $D$ of $(M^n, g, \phi)$ satisfies
\[    \lambda^2   \  \geq \  \frac{n}{4(n-2)}  \, \inf_M  S.      \]
 Equality occurs if and only if there exists a non-trivial spinor field  \, $\psi^{\ast}$ such that the differential equation
\begin{equation}   \nabla_X   \psi^{\ast}  =  - \frac{\lambda^{\ast}}{n} \,  X  \cdot  \psi^{\ast}  -   \frac{1}{n} \, \phi(X)  \cdot  \widetilde{D}  \psi^{\ast}   \end{equation}
holds for some real number $\lambda^{\ast}  \not=0  \in {\mathbb
R}$ and for all vector fields $X$.
\end{thm}

\bigskip  \noindent
\begin{thm}
Let  $( M^n, g, \phi), \, n \geq 3,$ be a locally decomposable
Riemannian spin manifold with positive scalar curvature $S > 0$.
Assume that $( M^n, g, \phi)$ is compact and
 \[ n_1 = {\rm dim} (T^+(M^n)) \ > \  n_2 = {\rm dim} (T^-(M^n)) \ \geq \ 1 .\]
Then any eigenvalue $\lambda$ of the Dirac operator $D$ of $(M^n,
g, \phi)$ satisfies
\[    \lambda^2   \  \geq \  \frac{n_1}{4(n_1-1)}  \, \inf_M  S.      \]
 Equality occurs if and only if there exists a non-trivial spinor field \, $\psi^{\ast}$ such that the differential equation
\begin{equation}   \nabla_X   \psi^{\ast}  =  - \frac{\lambda^{\ast}}{2 n_1} \,  X  \cdot  \psi^{\ast}  -   \frac{\lambda^{\ast}}{2  n_1} \, \phi(X)  \cdot  \psi^{\ast}   \end{equation}
holds for some real number $\lambda^{\ast} \not=0 \in {\mathbb R}$
and for all vector fields $X$.
\end{thm}

\bigskip  \noindent
Riemannian spin  manifolds admitting non-trivial solutions to the
field equation (1.3) (resp. the field equation (1.4)) are called
{\it limiting (locally decomposable) manifolds}.
In Section 4 we provide the following examples of limiting manifolds :     \\
Let $(M^{n_1}_1, g_1)$ and $(M^{n_2}_2, g_2)$ be Riemannian spin manifolds admitting Killing spinors. Let $( \overline{M}^{\, n_2}_2,  \overline{g}_2)$ be a Riemannian spin manifold  admitting parallel spinors.   For the classification of manifolds with  Killing spinors (resp. parallel spinors), we refer to [2, 22].      \\
(i)  If $n_1 = n_2  \geq 2$, then the Riemannian product manifold $(M^{n_1}_1 \times M^{n_2}_2, \, g_1 + g_2 )$ as well  as  $(M^{n_1}_1 \times  \overline{M}^{\,  n_2}_2, \, g_1 +  \overline{g}_2 )$   satisfies the limiting case of Theorem 1.1  (see Theorem 4.1 and 4.3).  \\
(ii)  If $n_1 > n_2  \geq 1$, then the Riemannian product manifold
$(M^{n_1}_1 \times  \overline{M}^{\, n_2}_2, \, g_1 +
\overline{g}_2 )$ satisfies the limiting case of Theorem 1.2 (see
Theorem 4.2 and 4.4).

\bigskip   \noindent
\section{Basic properties of the twisted Dirac operator $\widetilde{D}$ and some remarks on the limiting case of Theorem 1.1-1.2}

\noindent In this section we define the twist $\widetilde{D}$ of
the Dirac operator $D$, introduced in the previous section, in a
general setting. We then establish some formulas needed to prove
Theorem 1.1-1.2  in the next section.  We will  in fact show that
some basic results in  K\"{a}hler spin geometry can be translated
into the forms to be appropriate for almost product Riemannian
manifolds.  Because of the similarities between almost  product
Riemannian structures and almost  Hermitian structures, we will
describe the formulas in a unified way so as to be valid for both
types of structures.

\bigskip
Let $(M^n, g)$ be an n-dimensional  Riemannian spin manifold. Let
$\phi$ be a (1,1)-tensor field on $(M^n,  g)$ such that $\phi^2 =
\sigma  I,   \,  \sigma = \pm 1$, and
\[        g (  \phi X,   \phi  Y) =  g (X,   Y)            \]
for all vector fields $X, Y$ (Here $I$ stands for the identity
map).  Since
\[      g ( \phi X,  Y)   =  \sigma   \, g (X,  \phi Y)  ,    \]
the endomorphism $\phi$ is skew-symmetric if $\sigma = -1$ and
symmetric if $\sigma = 1$, respectively. Note that $(M^n, g,
\phi)$ is called an  {\it almost Hermitian manifold} if $\sigma =
-1$  and an  {\it almost product Riemannian manifold} if $\sigma =
1$, respectively. Let $\Sigma(M^n)$ be the spinor bundle of $(M^n,
g, \phi)$.  In terms of local  orthonomal frame field $(E_1,
\ldots, E_n)$, the spin derivative $\nabla $ and the Dirac
operator $D$, acting on sections $\psi  \in  \Gamma(\Sigma(M^n))$
of $\Sigma (M^n)$,
 are locally expressed as
\[             \nabla_X   \psi  = X(\psi) + \frac{1}{4} \sum_{i=1}^n
E_i \cdot \nabla_X  E_i \cdot \psi        \] and
\[     D \psi = \sum_{i=1}^n E_i \cdot \nabla_{E_i} \psi  ,   \]
respectively. Associated with the endomorphism $\phi$,  we define
the {\it  $\phi$-twist} of the Dirac operator $D$ by
\begin{equation}
\widetilde{D} \psi =  \sum_{i=1}^n E_i \cdot \nabla_{\phi(E_i)}
\psi   =  \sigma \, \sum_{i=1}^n \phi(E_i) \cdot \nabla_{E_i} \psi
.
\end{equation}
 Let $( \cdot, \cdot ) =  {\rm Re} \langle  \cdot, \cdot \rangle$ be the real part of the standard Hermitian product $\langle  \cdot, \cdot \rangle$ on the spinor bundle $\Sigma(M^n)$.
Then, from the identity
\begin{eqnarray}
&        &  -  \sigma \,   {\rm div}  \Big(  \sum_{i=1}^n   ( \psi_1,  \,   \phi(E_i)  \cdot  \psi_2 ) E_i  \Big)    \nonumber   \\
&        &     \nonumber    \\
&   =  &  (   \widetilde{D} \psi_1,  \, \psi_2 ) - ( \psi_1,  \,
\widetilde{D} \psi_2 )  - \sigma  \sum_{i=1}^n     ( \psi_1,  \,
( \nabla_{E_i}  \phi )(E_i)  \cdot   \psi_2  ),
\end{eqnarray}
we see that $\widetilde{D}$ is self-adjoint  with respect to
$L^2$-product  if $\phi$ satisfies
\[
{\rm div}(\phi)  = 0 .
\]

\noindent
\begin{pro}   On almost product Riemannian (resp. almost Hermitian) spin manifold $(M^n, g, \phi)$, we have
\begin{eqnarray*}
&       &   \widetilde{D}^2   \psi   -   D^2  \psi         \\
&       &         \\
&   =  &  -  \frac{1}{8}  \sum_{i,j,k,l=1}^n   g ( ( \nabla^2
\phi )(E_j, E_l, E_k) - ( \nabla^2  \phi )(E_j, E_k, E_l),  \,
\phi(E_i))
E_i  \cdot  E_j  \cdot  E_k  \cdot E_l   \cdot   \psi         \\
&      &       \\
&      &    +   \sum_{i,j=1}^n    \phi(E_i)    \cdot   (
\nabla_{E_i}  \phi ) (E_j)   \cdot   \nabla_{E_j}   \psi  ,
\end{eqnarray*}
where the second covariant derivative  $( \nabla^2  \phi )(Z, Y,
X)$ is defined by
\[   ( \nabla^2  \phi )(Z, Y, X)  =  \nabla_X  \{ ( \nabla_Y \phi)(Z)  \}   - (\nabla_Y \phi)(\nabla_X Z) - ( \nabla_{\nabla_X Y} \phi )(Z)  .  \]
\end{pro}

\noindent {\bf  Proof.}  Using a local orthonormal frame field
$(E_1, \ldots, E_n)$, we compute
\begin{eqnarray}
&     &   \widetilde{D}^2  \psi  = \sum_{i,j=1}^n E_i \cdot \nabla_{\phi E_i} ( E_j \cdot \nabla_{\phi E_j} \psi )    \nonumber   \\
&     &    \nonumber    \\
&  =  &  \sum_{i,j=1}^n  E_i \cdot E_j  \cdot  \nabla_{\phi E_i}
\nabla_{\phi E_j} \psi  + \sum_{i,j=1}^n  E_i  \cdot  \nabla_{\phi
E_i} E_j
\cdot  \nabla_{\phi E_j}  \psi    \nonumber   \\
&      &   \nonumber   \\
&   =  &  \frac{1}{2}  \sum_{i,j=1}^n    E_i \cdot E_j  \cdot  \Big\{ \nabla_{\phi E_i} \nabla_{\phi E_j} \psi  - \nabla_{\phi E_j} \nabla_{\phi E_i} \psi  \Big\}    \nonumber    \\
&       &   \nonumber    \\
&       &  -  \sum_{i=1}^n  \nabla_{\phi E_i} \nabla_{\phi E_i}
\psi   +  \sum_{i,j=1}^n  E_i  \cdot  \nabla_{\phi E_i} E_j
\cdot  \nabla_{\phi E_j}  \psi    \nonumber   \\
&       &    \nonumber   \\
&    = &   \frac{1}{2}  \sum_{i,j=1}^n  E_i  \cdot  E_j  \cdot  R
( \phi E_i, \phi  E_j )(\psi)   +  \triangle (\psi)  +
\frac{1}{2}    \sum_{i,j=1}^n  E_i  \cdot  E_j  \cdot
\nabla_{[ \phi E_i, \, \phi E_j ]}   \psi     \nonumber    \\
&       &     \nonumber    \\
&       &   -  \sum_{i=1}^n  \nabla_{ \nabla_{\phi E_i} ( \phi E_i
) }  \psi
 +  \sum_{i,j=1}^n  E_i  \cdot  \nabla_{\phi E_i} E_j
\cdot  \nabla_{\phi E_j}  \psi .
\end{eqnarray}
Since
\begin{eqnarray*}
&      &   \frac{1}{2}    \sum_{i,j=1}^n  E_i  \cdot  E_j  \cdot   \nabla_{[ \phi E_i, \, \phi E_j ]}   \psi      \\
&       &       \\
&  =  &       \sum_{i,j=1}^n  E_i  \cdot  E_j  \cdot
\nabla_{\nabla_{\phi E_i}  ( \phi E_j ) } \psi   +
 \sum_{i=1}^n  \nabla_{\nabla_{\phi E_i } ( \phi E_i ) } \psi      \\
&        &        \\
&    =    &     \sum_{i,j=1}^n  E_i  \cdot  E_j  \cdot
\nabla_{(\nabla_{\phi E_i}  \phi )  ( E_j ) } \psi    +
\sum_{i,j=1}^n  E_i  \cdot  E_j  \cdot   \nabla_{ \phi(
\nabla_{\phi E_i}  E_j )   } \psi  +
\sum_{i=1}^n  \nabla_{\nabla_{\phi E_i } ( \phi E_i ) } \psi      \\
&      &      \\
&   =   &   \sigma   \sum_{i,j=1}^n  E_i \cdot  ( \nabla_{\phi
E_i} \phi ) (E_j)   \cdot  \nabla_{E_j}  \psi   - \sum_{i,j=1}^n
E_i  \cdot  \nabla_{\phi E_i} E_j  \cdot  \nabla_{\phi E_j}  \psi
+  \sum_{i=1}^n  \nabla_{\nabla_{\phi E_i } ( \phi E_i ) } \psi
\end{eqnarray*}
and
\[    \sigma  \sum_{i,j=1}^n  E_i \cdot  ( \nabla_{\phi E_i} \phi ) (E_j)   \cdot  \nabla_{E_j}  \psi   =
\sum_{i,j=1}^n  \phi(E_i) \cdot  ( \nabla_{E_i} \phi ) (E_j)
\cdot  \nabla_{E_j}  \psi ,    \] the equation (2.3) now becomes
\begin{equation}
\widetilde{D}^2  \psi     =     \frac{1}{2}  \sum_{i,j=1}^n  E_i
\cdot  E_j  \cdot  R ( \phi E_i, \phi  E_j )(\psi)   +  \triangle
(\psi)    \nonumber
 +  \sum_{i,j=1}^n  \phi(E_i) \cdot  ( \nabla_{E_i} \phi ) (E_j)   \cdot  \nabla_{E_j}  \psi .
\end{equation}
On the other hand,  the identity
\[  ( \nabla^2 \phi ) (Z, Y, X) -   ( \nabla^2 \phi ) (Z, X, Y)  =  R(X, Y)( \phi Z)  - \phi \{ R(X, Y)(Z) \}     \]
implies
\begin{eqnarray*}
&      &   R( \phi E_i, \phi E_j, E_k, E_l ) -  R( E_i,  E_j, E_k, E_l )      \\
&  =   &  g(   ( \nabla^2 \phi ) (E_j,  E_l, E_k )  -   ( \nabla^2
\phi ) (E_j,  E_k, E_l ),   \  \phi E_i )
\end{eqnarray*}
and so
\begin{eqnarray}
&       &   \frac{1}{2}  \sum_{i,j=1}^n  E_i  \cdot  E_j  \cdot  R ( \phi E_i, \phi  E_j )(\psi)    \nonumber  \\
&       &    \nonumber    \\
&   =  &  - \frac{1}{8}  \sum_{i,j,k,l=1}^n    R( \phi E_i, \phi E_j, E_k, E_l )  E_i  \cdot E_j \cdot E_k \cdot E_l \cdot \psi   \nonumber   \\
&      &   \nonumber   \\
&    =  & - \frac{1}{8}  \sum_{i,j,k,l=1}^n    R( E_i,  E_j, E_k, E_l )  E_i  \cdot E_j \cdot E_k \cdot E_l \cdot \psi    \nonumber  \\
&      &   \nonumber   \\
&      &  - \frac{1}{8}  \sum_{i,j, k,l=1}^n  g (   ( \nabla^2
\phi ) (E_j,  E_l, E_k )  -   ( \nabla^2 \phi ) (E_j,  E_k, E_l ),
\  \phi E_i )
E_i \cdot E_j \cdot E_k \cdot E_l \cdot \psi   \nonumber  \\
&       &   \nonumber   \\
&   =  &  \frac{1}{4} \, S \, \psi      \\
&       &    \nonumber    \\
&       &   - \frac{1}{8}  \sum_{i,j, k,l=1}^n  g (   ( \nabla^2
\phi ) (E_j,  E_l, E_k )  -   ( \nabla^2 \phi ) (E_j,  E_k, E_l ),
\  \phi E_i )   E_i \cdot E_j \cdot E_k \cdot E_l \cdot \psi .
\nonumber
\end{eqnarray}
Applying (2.5) and the Schr\"{o}dinger-Lichnerowicz formula \[ D^2
= \triangle + \frac{1}{4} S \] to (2.4), we obtain the formula of
the proposition.      \hfill{QED.}

\noindent
\begin{pro}      On almost product Riemannian (resp. almost Hermitian) spin manifold $(M^n, g, \phi)$, we have
\begin{eqnarray*}
&       &    D ( \widetilde{D}  \psi )    +   \widetilde{D}  ( D \psi )          \\
&       &        \\
&   =  &  -  \sum_{i=1}^n    ( \nabla^2  \psi )  ( \phi E_i , E_i )  -  \sum_{i=1}^n    ( \nabla^2  \psi )  ( E_i ,   \phi E_i )      \\
&       &        \\
&       &  +   \frac{1}{2}   \sum_{i,j=1}^n    E_i   \cdot  E_j
\cdot   R (E_i,  \phi  E_j )(\psi)   +
 \frac{1}{2}   \sum_{i,j=1}^n     E_i   \cdot  E_j   \cdot   R ( \phi E_i,  E_j )(\psi)      \\
&       &       \\
&       &   +  \sigma   \sum_{i,j=1}^n     E_i    \cdot   (
\nabla_{E_i}  \phi ) (E_j)   \cdot   \nabla_{E_j}   \psi  ,
\end{eqnarray*}
where the second spinor derivative $( \nabla^2  \psi )(Y, X)$ is
defined by
\[    ( \nabla^2  \psi )(Y, X) = \nabla_X \nabla_Y  \psi -   \nabla_{\nabla_X Y} \psi .    \]
\end{pro}

\noindent {\bf Proof.}  Using a local orthonormal frame field
$(E_1, \ldots, E_n)$, we compute
\begin{eqnarray*}
&      &   D ( \widetilde{D}  \psi )    +   \widetilde{D}  ( D \psi )          \\
&       &        \\
&   =  &   \sum_{i,j=1}^n  E_i  \cdot  \nabla_{E_i} ( E_j \cdot
\nabla_{ \phi E_j } \psi )  +
 \sum_{i,j=1}^n  E_i  \cdot  \nabla_{\phi E_i} ( E_j \cdot  \nabla_{ E_j } \psi )     \\
&      &      \\
&  =  &  \sum_{i,j=1}^n  E_i \cdot E_j \cdot  \nabla_{E_i}
\nabla_{\phi E_j}  \psi  +
\sum_{i,j=1}^n  E_i \cdot E_j \cdot  \nabla_{\phi E_i} \nabla_{E_j}  \psi       \\
&      &      \\
&       &  + \sum_{i,j=1}^n  E_i \cdot \nabla_{E_i} E_j \cdot  \nabla_{\phi E_j}  \psi  +  \sum_{i,j=1}^n  E_i \cdot \nabla_{\phi E_i} E_j \cdot  \nabla_{E_j}  \psi        \\
&        &      \\
&   =   &  \frac{1}{2}   \sum_{i,j=1}^n  E_i \cdot E_j \cdot
\nabla_{E_i} \nabla_{\phi E_j}  \psi   -  \frac{1}{2}
\sum_{i,j=1}^n  E_j \cdot E_i
 \cdot  \nabla_{E_i} \nabla_{\phi E_j}  \psi      \\
&       &      \\
&        &  + \frac{1}{2} \sum_{i,j=1}^n  E_i \cdot E_j \cdot  \nabla_{\phi E_i} \nabla_{E_j}  \psi   -  \frac{1}{2} \sum_{i,j=1}^n  E_j \cdot E_i \cdot  \nabla_{\phi E_i} \nabla_{E_j}  \psi       \\
&         &      \\
&          &  -  \sum_{i=1}^n  \nabla_{E_i}  \nabla_{\phi E_i} \psi   -  \sum_{i=1}^n  \nabla_{\phi E_i}  \nabla_{E_i} \psi     \\
&      &      \\
&       &  + \sum_{i,j=1}^n  E_i \cdot \nabla_{E_i} E_j \cdot  \nabla_{\phi E_j}  \psi  +  \sum_{i,j=1}^n  E_i \cdot \nabla_{\phi E_i} E_j \cdot  \nabla_{E_j}  \psi        \\
&       &       \\
&   =   &  \frac{1}{2}  \sum_{i,j=1}^n  E_i \cdot E_j \cdot R( E_i,  \phi E_j)(\psi)  + \frac{1}{2}  \sum_{i,j=1}^n  E_i \cdot E_j \cdot R( \phi E_i,  E_j)(\psi)       \\
&        &       \\
&        &   - \sum_{i=1}^n  ( \nabla^2  \psi )( \phi E_i, E_i )   - \sum_{i=1}^n  ( \nabla^2  \psi )( E_i,  \phi E_i )       \\
&        &       \\
&        &  + \frac{1}{2}  \sum_{i,j=1}^n  E_i \cdot E_j \cdot
\nabla_{[ E_i, \, \phi E_j ]} \psi  +
\frac{1}{2}  \sum_{i,j=1}^n  E_i \cdot E_j \cdot  \nabla_{[ \phi E_i, \, E_j ]} \psi      \\
&         &      \\
&         & - \sum_{i=1}^n  \nabla_{ \nabla_{E_i} ( \phi E_i ) } \psi  -  \sum_{i=1}^n  \nabla_{ \nabla_{\phi E_i} E_i  } \psi     \\
&         &       \\
&       &  + \sum_{i,j=1}^n  E_i \cdot \nabla_{E_i} E_j \cdot  \nabla_{\phi E_j}  \psi  +  \sum_{i,j=1}^n  E_i \cdot \nabla_{\phi E_i} E_j \cdot  \nabla_{E_j}  \psi .       \\
\end{eqnarray*}
Since
\begin{eqnarray*}
&        &  \frac{1}{2}  \sum_{i,j=1}^n  E_i \cdot E_j \cdot
\nabla_{[ E_i, \, \phi E_j ]} \psi  +
\frac{1}{2}  \sum_{i,j=1}^n  E_i \cdot E_j \cdot  \nabla_{[ \phi E_i, \, E_j ]} \psi       \\
&        &      \\
&  =   &   \sum_{i,j=1}^n  E_i \cdot E_j \cdot  \nabla_{
\nabla_{E_i} ( \phi E_j ) } \psi   +
 \sum_{i,j=1}^n  E_i \cdot E_j \cdot  \nabla_{ \nabla_{\phi E_i}  E_j  } \psi        \\
&        &       \\
&        &   +  \sum_{i=1}^n  \nabla_{ \nabla_{E_i} ( \phi E_i ) } \psi  +  \sum_{i=1}^n  \nabla_{ \nabla_{\phi E_i} E_i  } \psi    \\
&        &      \\
&    =  &  \sum_{i,j=1}^n  E_i \cdot E_j \cdot  \nabla_{ ( \nabla_{E_i}  \phi )(E_j) } \psi   +  \sum_{i,j=1}^n  E_i \cdot E_j \cdot  \nabla_{ \phi ( \nabla_{E_i} E_j )} \psi      \\
&        &      \\
&        &  + \sum_{i,j=1}^n  E_i \cdot E_j \cdot  \nabla_{ \nabla_{\phi E_i}  E_j  } \psi   +  \sum_{i=1}^n  \nabla_{ \nabla_{E_i} ( \phi E_i ) } \psi  +  \sum_{i=1}^n  \nabla_{ \nabla_{\phi E_i} E_i  } \psi     \\
&        &      \\
&    =    &   \sigma \sum_{i,j=1}^n  E_i  \cdot  ( \nabla_{E_i} \phi )(E_j)  \cdot  \nabla_{E_j}  \psi     \\
&         &       \\
&         &   -  \sum_{i,j=1}^n  E_i \cdot \nabla_{E_i} E_j \cdot  \nabla_{\phi E_j}  \psi  -  \sum_{i,j=1}^n  E_i \cdot \nabla_{\phi E_i} E_j \cdot  \nabla_{E_j}  \psi        \\
&         &         \\
&         &   +  \sum_{i=1}^n  \nabla_{ \nabla_{E_i} ( \phi E_i )
} \psi  +  \sum_{i=1}^n  \nabla_{ \nabla_{\phi E_i} E_i  } \psi,
\end{eqnarray*}
we now obtain the formula of the proposition. \hfill{QED.}

\bigskip  \noindent
{\bf Remark.}   On K\"{a}hler spin manifolds (i.e., if  $\phi$ is
skew-symmetric and   $\nabla \phi = 0$),  the relations in
Proposition 2.1-2.2 simplify to the well-known relation [13]
\begin{equation}
\widetilde{D}^2 = D^2
\end{equation}
  and
\begin{equation}      D \widetilde{D} +  \widetilde{D} D = 0 ,
\end{equation}
respectively. However, the relation (2.7) does not generally hold
on locally decomposable Riemannian spin manifolds.

\bigskip
Analogously to the {\it  K\"{a}hlerian twistor equation} [14, 15],
we now  consider the following spinor field equation
\[
\nabla_{\phi X}  \psi  = p \, \phi(X)  \cdot  D \psi  +  q \, X
\cdot  \widetilde{D}  \psi ,    \qquad  p \not= 0,  \,  q \not= 0
\in  {\mathbb R},
\]
which is equivalent to
\begin{equation}
\nabla_X   \psi  = p \,  X  \cdot  D \psi  +  \sigma  \,  q \,
\phi(X)  \cdot  \widetilde{D}  \psi.
\end{equation}

\bigskip  \noindent
{\bf  Definition 2.1}  A non-trivial solution $\psi$ to the field
equation (2.8) on almost product Riemannian (resp. almost
Hermitian) spin manifold $(M^n, g, \phi)$ is called {\it
quasi-twistor-spinor} (resp. {\it Hermitian twistor-spinor}) of
type (p,q).

\noindent
\begin{pro}
Let $(M^n,  g, \phi)$ admit a  quasi-twistor-spinor $\psi$ (resp.
Hermitian twistor-spinor) of type (p, q).
  Then we have
\begin{eqnarray}
\frac{1}{2}  {\rm  Ric}(X)  \cdot  \psi   & = & - p \, X  \cdot
D^2  \psi - (2p+1) \nabla_X (D \psi)  - 2 \sigma q  \nabla_{\phi
X}
( \widetilde{D} \psi ) - \sigma  q \, \phi(X)  \cdot  D  \widetilde{D} \psi         \nonumber      \\
&      &      \nonumber       \\
&      &     +  \sigma  q  \sum_{i=1}^n   E_i  \cdot  (
\nabla_{E_i} \phi )(X)  \cdot   \widetilde{D}  \psi .
\end{eqnarray}
Contracting this equation via $\displaystyle  S \, \psi = -
\sum_{i=1}^n E_i  \cdot  {\rm Ric}(E_i)  \cdot  \psi$ gives
\begin{eqnarray}
\frac{1}{2} \, S \, \psi  & = & (1+2p-np) D^2 \psi  + 2 \sigma q
\,   \widetilde{D}^2   \psi   +  \sigma  q  \sum_{i=1}^n
  E_i  \cdot \phi(E_i)  \cdot  D  \widetilde{D}  \psi            \nonumber     \\
&         &        \nonumber     \\
&         &    -  \sigma  q  \sum_{i,j=1}^n   E_j  \cdot  E_i
\cdot ( \nabla_{E_i}  \phi ) (E_j)   \cdot  \widetilde{D}  \psi .
\end{eqnarray}
\end{pro}

\noindent {\bf Proof.}  Applying (2.8) to the $(\frac{1}{2}
Ricci)$-formula (see Lemma 1.2 in [10])
\[    \frac{1}{2} {\rm Ric} (X) \cdot \psi = D( \nabla_X \psi)  -  \nabla_X (D \psi) -  \sum_{i=1}^n E_i  \cdot  \nabla_{\nabla_{E_i} X} \psi ,  \]
we have
\begin{eqnarray*}
   \frac{1}{2}  {\rm Ric}(X)   \cdot   \psi
&  =  & -  \nabla_X (D \psi)  +  D \Big( p X \cdot D\psi  + \sigma q \, \phi(X)  \cdot  \widetilde{D} \psi  \Big)      \\
&       &      \\
&       & -   \sum_{i=1}^n  E_i  \cdot  \Big\{  p \nabla_{E_i} X \cdot D \psi + \sigma q \, \phi( \nabla_{E_i} X )  \cdot  \widetilde{D} \psi  \Big\}       \\
&        &        \\
&   =   &  -  \nabla_X (D \psi)   +  p \sum_{i=1}^n  E_i  \cdot \nabla_{E_i} X \cdot D \psi - 2p \nabla_X ( D \psi ) - p X \cdot D^2 \psi      \\
&        &      \\
&        &  + \sigma q  \sum_{i=1}^n  E_i  \cdot  \nabla_{E_i}
(\phi X)  \cdot \widetilde{D} \psi - 2 \sigma q \nabla_{\phi X} (
\widetilde{D} \psi )
- \sigma q \, \phi(X) \cdot  D \widetilde{D} \psi      \\
&        &      \\
&       &     -  p \sum_{i=1}^n  E_i  \cdot  \nabla_{E_i} X \cdot D \psi - \sigma q  \sum_{i=1}^n E_i \cdot \phi( \nabla_{E_i} X )  \cdot  \widetilde{D} \psi         \\
&       &      \\
&   =   &  - p \, X  \cdot  D^2  \psi - (2p+1) \nabla_X (D \psi)
- 2 \sigma q  \nabla_{\phi X}
( \widetilde{D} \psi ) - \sigma  q \, \phi(X)  \cdot  D  \widetilde{D} \psi         \nonumber      \\
&      &      \nonumber       \\
&      &     +  \sigma  q  \sum_{i=1}^n   E_i  \cdot  (
\nabla_{E_i} \phi )(X)  \cdot   \widetilde{D}  \psi .
\end{eqnarray*}
\hfill{QED.}

\noindent
\begin{cor}  Let $(M^n, g, \phi)$ be an almost product Riemannian (resp. almost Hermitian) spin manifold with $\nabla \phi =0$. Assume that
$(M^n,  g, \phi)$ admits a  quasi-twistor-spinor $\psi$ (resp.
Hermitian twistor-spinor) of type (p, q).  Then we have
\begin{equation} 4(p+q+1)  \,   D^2  \psi  =  S \,  \psi .  \end{equation}
\end{cor}

\noindent {\bf  Proof.}  We first consider the case $\sigma =-1$.
Let $\Omega$ be the  fundamental 2-form defined by
\[   \Omega(X, Y) = g( X, \phi Y)  .    \]
Then we know that
\begin{equation}
D \Omega - \Omega D = - 2 \widetilde{D},    \qquad
\widetilde{D} \Omega -  \Omega \widetilde{D} =  2 D .
\end{equation}
On the other hand, contracting the equation (2.8), we obtain
\begin{equation}
(1+nq) \, \widetilde{D} \psi = - 2p \, \Omega \cdot D \psi ,
\qquad   (1+np) \, D \psi = 2 q \, \Omega \cdot  \widetilde{D}
\psi .
\end{equation}
Applying (2.12) to (2.13), we have
\begin{equation}
\Omega \cdot D \widetilde{D} \psi =  \frac{1+4p +nq}{2p} \, D^2
\psi = \frac{1+np+4q}{2q} \, D^2 \psi .
\end{equation}
Consequently, inserting (2.14) into (2.10) gives the formula
(2.11) of the corollary. Now we prove that (2.11) is also true for
the other case $\sigma =1$. Contracting the equation (2.8), we
obtain
\begin{equation}
(1+np) \, D \psi = - q  \,  {\rm Tr}(\phi) \, \widetilde{D} \psi
,   \qquad   (1+nq) \, \widetilde{D} \psi = - p  \,  {\rm
Tr}(\phi) \, D \psi  .
\end{equation}
Inserting (2.15) into (2.10) gives the formula (2.11).
\hfill{QED.}

\bigskip
In the rest of the section we make some remarks on the limiting
case of Theorem 1.1 (resp. Theorem 1.2).  It is obvious that a
locally decomposable Riemannian spin manifold $(M^n, g, \phi)$ is
a limiting manifold of Theorem 1.1 if and only if $(M^n, g, \phi)$
admits such an eigenspinor $\psi^{\ast}$ of the Dirac operator
that is a  quasi-twistor-spinor of type $(- \frac{1}{n}, -
\frac{1}{n})$.   Therefore, by Corollary 2.1,  the scalar
curvature of any limiting manifold of Theorem 1.1 is necessarily
constant.

\bigskip
In order to discuss the limiting case of Theorem 1.2, we now
consider a special type of spinor field equation
\begin{equation}
\nabla_X \psi = a \, X \cdot  \psi  +  b \, \phi(X)  \cdot  \psi ,
\qquad  a \not=0, \, b \not= 0  \in  {\mathbb R},
\end{equation}
which is closely related to the quasi-twistor equation  (2.8).
From (2.15) we observe that  the equation (2.8) reduces to (2.16) in the following cases :   \\
(i) $D \psi = \lambda  \psi$  for some $\lambda \not=0  \in {\mathbb R}$ and ${\rm Tr}(\phi)  \not=0$.     \\
(ii) $D \psi = \lambda  \psi$  for some $\lambda \not=0  \in {\mathbb R}$ and $q  \not= - \frac{1}{n}$.     \\

\bigskip \noindent
{\bf Definition 2.2} A non-trivial solution $\psi$ to the field
equation (2.16) on almost product Riemannian spin manifold $(M^n,
g, \phi)$ is called {\it quasi-Killing spinor} of type (a,b).

\bigskip  \noindent
\noindent
\begin{pro}
Let $(M^n,  g, \phi)$ be an almost product Riemannian spin
manifold admitting a  quasi-Killing spinor $\psi$  of type (a,b).
  Then we have
\begin{eqnarray}
&       &   {\rm  Ric}(X)  \cdot  \psi   =    4  \Big\{ (n-1) a^2  + ab \cdot {\rm Tr}(\phi) - b^2  \Big\}  X \cdot  \psi      \\
&      &   \nonumber    \\
&      &   \qquad  \quad + 4 \Big\{ b^2  \cdot {\rm Tr}(\phi)  +
(n-2) ab  \Big\} \phi(X)  \cdot    \psi   + 2 b  \sum_{i=1}^n  E_i
\cdot  ( \nabla_{E_i} \phi )(X)  \cdot  \psi .   \nonumber
\end{eqnarray}
In particular, the scalar curvature $S$ must be constant and given
by
\begin{equation}
S =  4n(n-1) a^2  + 8(n-1) ab \cdot {\rm Tr}(\phi)  - 4n b^2 +  4
b^2 ( {\rm Tr} \phi )^2 .
\end{equation}
\end{pro}

\noindent {\bf Proof.}  Applying (2.16) to the $(\frac{1}{2}
Ricci)$-formula, we obtain the equation (2.17) immediately.
Contracting (2.17) gives
\begin{eqnarray*}
 S \, \psi & = & 4n  \Big\{ (n-1) a^2 + ab \cdot {\rm Tr}(\phi) - b^2  \Big\}  \psi    \\
&      &      \\
&      & + 4  \, {\rm Tr}(\phi)  \Big\{ b^2 \cdot {\rm Tr}(\phi) +
(n-2) ab  \Big\}  \psi + 4b \, {\rm div}(\phi) \cdot  \psi .
\end{eqnarray*}
Thus ${\rm div}(\phi)  =0$ must vanish, and rewriting gives
(2.18).     \hfill{QED.}

\bigskip  \noindent
Certainly, a locally decomposable Riemannian spin manifold $(M^n,
g, \phi)$ is a limiting manifold of Theorem 1.2  if and only if
$(M^n, g, \phi)$ admits a  quasi-Killing spinor $\psi^{\ast}$ of
type $(- \frac{\lambda^{\ast}}{2 n_1}, - \frac{\lambda^{\ast}}{2
n_1}), \, \lambda^{\ast}  \not=0  \in {\mathbb R}$. In this case,
the Ricci tensor (2.17) and the scalar curvature (2.18)  simplify
to
\begin{equation}
{\rm Ric} (X, Y)   =  \frac{2(n_1-1)}{n_1^2}  \cdot
(\lambda^{\ast})^2 \cdot \Big\{ g(X, Y)  + g( \phi(X), Y)  \Big\}
\end{equation} and
\[   S =  \frac{4(n_1-1)}{n_1} \cdot  (\lambda^{\ast})^2  , \]
respectively. If $(M^n, g) = (M_1^{n_1}  \times  M_2^{n_2},  \,
g_1  +  g_2 )$ is a global Riemannian product, then (2.19) implies
that $(M_1^{n_1}, g_1)$ is necessarily Einstein with positive
scalar curvature and $(M_2^{n_2}, g_2 )$ is Ricci-flat.

\bigskip   \noindent
\section{Proof of Theorem 1.1-1.2}

\noindent Let  $( M^n, g, \phi)$ be a locally decomposable
Riemannian spin manifold. Since $\nabla \phi = 0$ vanishes on
$M^n$, the operator $\widetilde{D}$ is self-adjoint with respect
to $L^2$-product (see (2.2)). Moreover,  Proposition 2.1 implies
\begin{equation}
\widetilde{D}^2 = D^2  .
\end{equation}
Let us define the {\it quasi-twistor operator}  \, $\mathcal{T} :
\Gamma ( T(M^n) )   \times  \Gamma( \Sigma (M^n) )
\longrightarrow \Gamma ( \Sigma(M^n))$ by
\[   \mathcal{T}_X  (\varphi)  : =  \nabla_X  \varphi - p \, X  \cdot  D \varphi  - q  \, \phi(X)  \cdot  \widetilde{D} \varphi ,  \qquad   p \not=0, \, q\not=0   \in  {\mathbb R}.     \]
Then a direct calculation using (3.1) yields
\begin{eqnarray}
0   &  \leq  &  \int_{M^n}  \sum_{i=1}^n  ( \mathcal{T}_{E_i}  (\varphi)   ,  \,  \mathcal{T}_{E_i}  (\varphi) ) \mu       \\
&     &   \nonumber     \\
&  =  &   \int_{M^n}  \Big\{  ( \nabla \varphi, \, \nabla \varphi)
+ (n p^2 + 2p + n q^2 + 2q ) ( D^2 \varphi, \,  \varphi )  + 2 pq
\cdot  {\rm Tr}(\phi) \cdot ( \widetilde{D} D \varphi, \,  \varphi
)    \Big\}  \mu,      \nonumber
\end{eqnarray}
where $\mu$ is the volume form of $(M^n, g)$. Now assume that
\[ n_1 =n_2,    \quad     \mbox{i.e.,}  \quad  {\rm Tr}(\phi) = 0.    \]
Then the equations in  (2.15) imply that the optimal parameters
$p, q$ are
\[  p = q = - \frac{1}{n}.   \]
Let  $\psi$  be an eigenspinor of $D$ with eigenvalue $\lambda$.
Then, applying the Schr\"{o}dinger-Lichnerowicz formula  to the
equation (3.2), we obtain
\[
0   \  \leq  \  \int_{M^n}  \sum_{i=1}^n  ( \mathcal{T}_{E_i}
(\psi)   ,  \,  \mathcal{T}_{E_i}  (\psi) ) \mu
     \  =  \   \int_{M^n}  \Big\{  \frac{n-2}{n} \cdot  \lambda^2  -  \frac{S}{4}  \Big\}(\psi, \psi)  \mu .
\]
This proves the inequality of Theorem 1.1. The limiting case of
the inequality is clear.

\bigskip
Next we prove Theorem 1.2. In order to control the last term $2 pq
\cdot  {\rm Tr}(\phi) \cdot ( \widetilde{D} D \varphi, \,  \varphi
) $ in (3.2), we introduce free parameters $a,b \in {\mathbb R}$
and compute
\begin{eqnarray}
0   &  \leq  &  \int_{M^n}  \Big\{ a^2 ( \, \widetilde{D} \varphi - b D \varphi,  \, \widetilde{D} \varphi - b D \varphi \, ) + \sum_{i=1}^n  ( \mathcal{T}_{E_i}  (\varphi)   ,  \,  \mathcal{T}_{E_i}  (\varphi) ) \Big\}  \mu     \nonumber   \\
&     &   \nonumber     \\
& =  &  \int_{M^n}   \Big\langle  - \frac{S}{4} ( \varphi,  \varphi )  + (1+ n p^2 + 2p + n q^2 + 2q + a^2 + a^2 b^2 ) ( D^2 \varphi, \,  \varphi )   \nonumber   \\
&     &     \nonumber    \\
&    &   \qquad  +    \{ 2 pq \,  {\rm Tr}(\phi) - 2 a^2 b \} (
\widetilde{D} D \varphi, \,  \varphi )    \Big\rangle  \mu,
\end{eqnarray}
Here we choose the parameters $a \not= 0, \, b \not= 0$ in such a
way that the last term in (3.3) vanishes and the equations in
(2.15) are satisfied with  $\widetilde{D} \psi  =  b D \psi $,
i.e.,
\begin{equation}     b =  \frac{ pq \, {\rm Tr}(\phi)}{a^2}  = -  \frac{ p \, {\rm Tr}(\phi) }{1+nq} =  -  \frac{1+np}{q \, {\rm Tr}(\phi)} .   \end{equation}
From this, it follows immediately that
\begin{equation}       ( {\rm Tr} \phi )^2  = (n_1 - n_2)^2 = \frac{(1+np)(1+nq)}{pq}     \end{equation}
and
\begin{equation}
a^2 =  - q(1+nq) ,    \qquad    b^2 =   \frac{ p (1+np)}{q(1+nq)}
.
\end{equation}
Since $a^2 > 0$ and  $ b^2  > 0$, we see that  $- \frac{1}{n} < p
< 0 $  and $- \frac{1}{n} < q < 0 .$ Inserting (3.6) into (3.3)
and assuming that $\varphi = \psi$ is an eigenspinor of $D$ with
eigenvalue $\lambda$, we now find that
\begin{eqnarray}
0   &  \leq  &  \int_{M^n}  \Big\{ a^2 ( \, \widetilde{D} \psi - b D \psi,  \, \widetilde{D} \psi - b D \psi \, ) + \sum_{i=1}^n  ( \mathcal{T}_{E_i}  (\psi)   ,  \,  \mathcal{T}_{E_i}  (\psi) ) \Big\}  \mu       \nonumber    \\
&     &   \nonumber     \\
& =  &  \int_{M^n}   \Big\{  (1+p+q) \lambda^2  -  \frac{S}{4}
\Big\} (\psi, \psi)  \mu .
\end{eqnarray}
Applying Lagrange's method to the function $f(p,q) := 1 +p+q$ with
the side condition (3.5), one verifies easily that $f(p,q)$ has
its minimum $\frac{n_1-1}{n_1}$ at the point
\begin{equation}     \Big( p= - \frac{1}{2  n_1}, \, q= - \frac{1}{2  n_1}  \Big)  .     \end{equation}
Consequently, (3.7) together with (3.8) leads us to the inequality
of Theorem 1.2. Moreover, inserting (3.8) into (3.6) gives
\begin{equation}
a^2 =  \frac{n_1- n_2}{4  \,  n_1^2},    \qquad   b=1 .
\end{equation}
Thus the limiting case of Theorem 1.2 is clear from (3.7).

\bigskip   \noindent
\section{Some limiting manifolds}

\noindent We show that some special types of Riemannian product
manifolds satisfy the limiting case of Theorem 1.1 (resp. Theorem
1.2). For that purpose we need to recall some algebraic formulas
describing  the action of  the Clifford algebra on tensor products
of spinor fields [4, 10, 17]. We begin with the case that  the
first manifold $(M_1^{2m} , g_1)$ is of even dimension $2m \geq 2$
and the second manifold $(M_2^{r} , g_2)$ is of  general dimension
$r \geq 2$.  In this case the Riemannian product manifold
$(M_1^{2m} \times M_2^{r} , \, g_1 + g_2)$  possesses  a naturally
induced spin structure and the spinor bundle of  $(M_1^{2m} \times
M_2^{r} , \, g_1 + g_2)$   is no other than the tensor product of
the spinor bundle of $(M_1^{2m}, g_1)$ and the spinor bundle of
$(M_2^{r}, g_2)$. Therefore, if $\psi_1$ and $\psi_2$ are a spinor
field on $(M_1^{2m} , g_1)$ and $(M_2^{r} , g_2)$, respectively,
then the tensor product $\psi_1 \otimes \psi_2$ is well defined on
$(M_1^{2m} \times M_2^{r} , \, g_1 + g_2)$. Let us denote by $(E_1
, \ldots , E_{2m})$ and $(F_1 , \ldots , F_{r})$ a local
orthonormal frame on $(M_1^{2m} , g_1)$ and $(M_2^{r} , g_2)$,
respectively. Identifying $(E_1 , \ldots , E_{2m})$ and $(F_1 ,
\ldots , F_{r})$ with their lifts to $(M_1^{2m} \times M_2^{r} ,
\, g_1 + g_2)$, we regard $(E_1 , \ldots , E_{2m} , F_1 , \ldots ,
F_{r})$ as a local orthonormal frame on  $(M_1^{2m} \times M_2^{r}
, g_1 + g_2)$. Then the Clifford bundle ${\rm Cl}(M_1^{2m}  \times
M_2^{r})$ of   $(M_1^{2m} \times M_2^{r} , g_1 + g_2)$ acts on the
spinor bundle
  $ \Sigma(M_1^{2m} \times  M_2^{r}) $ via
\begin{eqnarray}
E_i \cdot ( \psi_1  \otimes  \psi_2 )  &  =  & ( E_i  \cdot  \psi_1 )  \otimes  \psi_2  ,   \qquad  1 \leq i \leq 2m  ,  \\
&      &        \nonumber      \\
F_j  \cdot  ( \psi_1  \otimes  \psi_2 )  &  =  &  {(\sqrt{-1})}^m
( \mu_1 \cdot \psi_1 )   \otimes  ( F_j \cdot \psi_2 ) ,   \qquad
1 \leq  j \leq r,
\end{eqnarray}
where  $\mu_1 = E^1 \wedge \cdots \wedge E^{2m}, \, E^k := g_1 (
E_k,  \, \cdot  \, ),$ is the volume form of $(M_1^{2m} , g_1)$.
Denote by $\nabla^1$ (resp. $ \nabla^2$)  the Levi-Civita
connection and by  $D_1 $ (resp. $D_2$) the Dirac operator of
$(M_1^{2m} , g_1)$ (resp. $(M_2^{r} , g _2)$).  From (4.1)-(4.2),
we  immediately obtain the following formulas for the spin
derivative $\nabla$ and the Dirac operator $D$ of $(M_1^{2m}
\times M_2^{r} , \,  g_1 + g_2)$ :
\begin{eqnarray}
\nabla_X (\psi_1 \otimes \psi_2) & = & (\nabla^1_{\pi_1 (X)} \,
\psi_1) \otimes \psi_2 +
\psi_1 \otimes (\nabla^2_{\pi_2 (X)} \, \psi_2 ) ,           \\
&    &     \nonumber    \\
D (\psi_1 \otimes \psi_2) & = & (D_1 \psi_1) \otimes \psi_2  +
{(\sqrt{-1})}^m
(\mu_1 \cdot \psi_1) \otimes (D_2 \psi_2)  ,      \\
&    &     \nonumber     \\
D^2 (\psi_1 \otimes \psi_2) & = & ( D_1^2 \psi_1 )   \otimes
\psi_2 +  \psi_1 \otimes ( D_2^2 \psi_2 )    ,
\end{eqnarray}
where \, $\pi_1 : T(M_1^{2m} \times M_2^{r}) \longrightarrow
T(M_1^{2m})$ and $\pi_2 : T(M_1^{2m} \times M_2^{r})
\longrightarrow T(M_2^{r})$ are the natural projections.

\noindent
\begin{thm}
Let $(M_1^{2m}, g_1), \, m \geq 1,$ be a Riemannian spin manifold
admitting a Killing spinor $\psi_1$ with $D_1 \psi_1 = \lambda_1
\psi_1, \, \lambda_1 \not=0  \in  {\mathbb R}$. Let $(M_2^{2m},
g_2)$ be a Riemannian spin manifold admitting a Killing spinor
$\psi_2$ with $D_2 \psi_2 = \lambda_2, \, \lambda_2  \in {\mathbb
R}$ (Here we allow $\lambda_2$ to be zero).  Then the Riemannian
product manifold $(M_1^{2m} \times M_2^{2m} , \, g_1 + g_2)$
admits such a non-trivial eigenspinor $\psi^{\ast}$ of the Dirac
operator $D$ that satisfies the equation (1.3) of Theorem 1.1.
\end{thm}

\noindent {\bf Proof.} The spinor bundle $\Sigma(M_1^{2m})$ of
$(M_1^{2m}, g_1)$ decomposes into $\Sigma(M_1^{2m}) = \Sigma^+
(M_1^{2m}) \oplus \Sigma^- (M_1^{2m})$ under the action of the
volume element $\mu_1 = E_1  \cdots  E_{2m}$ ,
\[ \Sigma^{\pm} (M_1^{2m}) := \{ \, \varphi \in \Sigma (M_1^{2m}) : \mu_1 \cdot \varphi = \pm (\sqrt{-1})^m \varphi  \, \} . \]
Let  $\psi_1^{\pm} \in \Gamma (\Sigma^{\pm} (M_1^{2m}))$ be the
positive and negative part of  $\psi_1$, respectively.   Set
\[    \psi : = \psi_1^+ \otimes \psi_2 .     \]
Then, using the formula (4.5), we have
\[    D^2 \psi  =  ( \lambda^2_1 + \lambda^2_2 ) (\psi^+_1  \otimes \psi_2 ) =  ( \lambda^2_1 + \lambda^2_2 )   \psi.   \]
Let $S_1$ and $S_2$ be the scalar curvature of $(M_1^{2m}, g_1)$
and $(M_2^{2m}, g_2)$, respectively. Then $\lambda^{\ast} : =
\sqrt{\lambda^2_1 + \lambda^2_2}$ is related to the scalar
curvature $S = S_1 + S_2$ of $(M_1^{2m} \times M_2^{2m} , \, g_1 +
g_2 )$ as
\[   (\lambda^{\ast})^2 = \frac{1}{4} \cdot   \frac{2m}{2m-1}  \cdot  S_1   +  \frac{1}{4} \cdot   \frac{2m}{2m-1}  \cdot  S_2   =
\frac{1}{4}  \cdot  \frac{n}{n-2}  \cdot  S ,   \quad  n= 4m , \]
and
\[
\psi^{\ast}  \,  :=  \,  \lambda^{\ast}  \, \psi + D \psi
 =    \{ \lambda^{\ast}  + \lambda_2 (-1)^m \} ( \psi_1^+ \otimes \psi_2 ) + \lambda_1 ( \psi_1^- \otimes \psi_2 )
\]
is indeed such a non-trivial eigenspinor of $D$ with eigenvalue
$\lambda^{\ast}$ that satisfies the equation (1.3). \hfill{QED.}

\noindent
\begin{thm}
Let $(M_1^{2m}, g_1), \, m \geq 1,$ be a Riemannian spin manifold
admitting a Killing spinor $\psi_1$ with $D_1 \psi_1 = \lambda_1
\psi_1, \, \lambda_1 \not=0  \in  {\mathbb R}$. Let $(
\overline{M}_2^{\, r},  \overline{g}_2), \, 2m > r \geq 2,$ be a
Riemannian spin manifold admitting a parallel spinor $\psi_2$. Let
$( S^1, g_S)$ be a circle with the standard metric. Then the
Riemannian product manifold $(M_1^{2m} \times \overline{M}_2^{\,
r} , \, g_1 +  \overline{g}_2)$ as well as $( M_1^{2m}  \times
S^1, \, g_1 + g_S )$  admits such a non-trivial eigenspinor
$\psi^{\ast}$ of the Dirac operator $D$ that satisfies the
equation (1.4) of Theorem 1.2.
\end{thm}

\noindent {\bf  Proof.}   Applying the argument  in the proof  for
Theorem 4.1,  one proves that $(M_1^{2m} \times \overline{M}_2^{\,
r} , \, g_1 +  \overline{g}_2)$ satisfies the limiting case of
Theorem 1.2.    To prove the latter part of the theorem,  it
suffices to check that the lift $\psi_1^{\ast}$ of  $\psi_1$ to $(
M_1^{2m}  \times S^1, \, g_1 + g_S )$  satisfies the equation
(1.4) of Theorem 1.2.   \hfill{QED.}

\bigskip
Now we proceed to the other case that the first manifold
$(M_1^{2m+1} , \, g_1),  \,  m  \geq 1$, as well as  the second
manifold $(M_2^{2s+1} , \, g_2),  \, s \geq 0$,  is of   odd
dimension. We will modify the relations (4.1)-(4.5)  slightly to
be appropriate in this case. Let $(N^1, g_N)$ be a 1-dimensional
connected manifold (i.e., an open interval or a circle) with the
standard metric $g_N$,  and  let $( Q^{2m+2} = M_1^{2m+1}  \times
N^1,  \,  g_Q = g_1  + g_N )$ be the  Riemannian product manifold.
Let $E_{2m+2}$ denote a unit vector field on $( N^1, g_N)$ as well
as the lift to $(Q^{2m+2}, g_Q)$.  Denote by $(E_1 , \ldots ,
E_{2m+1})$ a local orthonormal frame on $(M_1^{2m+1} , g_1)$ as
well as the lift to $(Q^{2m+2}, g_Q)$. Then the spinor  bundle
$\Sigma(M_1^{2m+1})$ of $(M_1^{2m+1} , \, g_1)$ may be thought to
be embedded in the positive part $\Sigma^+( Q^{2m+2} )$ (resp.  in
the negative part $\Sigma^- (Q^{2m+2})$ ) of the spinor bundle
$\Sigma( Q^{2m+2} )$ of $( Q^{2m+2}, g_Q)$, the Clifford
multiplication ${\rm Cl}( M_1^{2m+1} )  \times \Sigma (M_1^{2m+1})
\longrightarrow  \Sigma (M_1^{2m+1})$ being
 naturally related to the one ${\rm Cl}( Q^{2m+2} )  \times \Sigma (Q^{2m+2})   \longrightarrow  \Sigma (Q^{2m+2})$ via
\begin{equation}    E_i  \cdot  ( \pi_Q \psi^{\pm} )  =  \pi_Q ( E_i \cdot  E_{2m+2} \cdot  \psi^{\pm}) ,    \qquad   1 \leq i  \leq 2m+1  ,    \end{equation}
where $\psi^{\pm}  \in  \Gamma(\Sigma^{\pm} (Q^{2m+2}))$ and
$\pi_Q : \Sigma^{\pm} (Q^{2m+2})  \longrightarrow
\Sigma(M_1^{2m+1})$ is the restriction map. Let $(F_1 , \ldots ,
F_{2s+1})$ be a local orthonormal frame on $(M_2^{2s+1} , g_2)$.
Identifying $(E_1 , \ldots , E_{2m+1})$ and $(F_1 , \ldots ,
F_{2s+1})$ with their lifts to $(M_1^{2m+1} \times M_2^{2s+1} , \,
g_1 + g_2)$, we  regard $(E_1 , \ldots , E_{2m+1} , F_1 , \ldots ,
F_{2s+1})$ as a local orthonormal frame on  $(M_1^{2m+1} \times
M_2^{2s+1} , \, g_1 + g_2)$. Then,  with help of (4.6),  one can
define a natural action of the Clifford bundle ${\rm Cl}
(M_1^{2m+1}  \times M_2^{2s+1})$ of  $( M_1^{2m+1}  \times
M_2^{2s+1},  \, g_1 + g_2 )$ on the spinor bundle
\begin{eqnarray*}
\Sigma( M_1^{2m+1}  \times M_2^{2s+1})  & =  &  \Big\{ \Sigma(M_1^{2m+1})  \oplus \Sigma (M_1^{2m+1})  \Big\}  \otimes \Sigma (M_2^{2s+1})  \\
&       &     \\
&  \subset   &  \Big\{ \Sigma^+ (Q^{2m+2})  \oplus  \Sigma^-
(Q^{2m+2})  \Big\}   \otimes \Sigma (M_2^{2s+1})
\end{eqnarray*}
of  $( M_1^{2m+1}  \times M_2^{2s+1},  \, g_1 + g_2)$  by
\begin{eqnarray}
E_i \cdot  \{ (\psi_1^+ + \psi_1^-)  \otimes  \psi_2 \}  &  =  & \{ E_i  \cdot  E_{2m+2} \cdot ( \psi_1^+  +  \psi_1^- )   \}  \otimes  \psi_2 ,     \quad    1 \leq i \leq 2m+1,     \\
&        &     \nonumber    \\
F_j  \cdot  \{  (  \psi_1^+  +  \psi_1^- )  \otimes  \psi_2  \}  &
=  &   \{  E_{2m+2}   \cdot  ( \psi_1^- -   \psi_1^+ )   \}
\otimes ( F_j  \cdot  \psi_2 )  ,    \quad   1 \leq  j  \leq
2s+1,
\end{eqnarray}
where $\psi_1^{\pm} \in \Gamma( \Sigma (M_1^{2m+1})) \subset
\Gamma( \Sigma^{\pm} (Q^{2m+2}))$ and  $ \psi_2  \in  \Gamma(
\Sigma( M_2^{2s+1})) $. Denote by $\nabla^1$ (resp. $ \nabla^2$)
the Levi-Civita connection and by  $D_1 $ (resp. $D_2$) the Dirac
operator of $(M_1^{2m+1} , g_1)$ (resp. $(M_2^{2s+1} , g _2)$).
From (4.7)-(4.8), we now obtain the following formulas for the
spin derivative $\nabla$ and the Dirac operator $D$ of
$(M_1^{2m+1} \times M_2^{2s+1} , \,  g_1 + g_2)$ :
\begin{eqnarray}
\nabla_X ( (\psi_1^+ + \psi_1^-)  \otimes \psi_2) & = & \Big\{
\nabla^1_{\pi_1 (X)} \, \psi_1^+ +  \nabla^1_{\pi_1 (X)} \psi_1^-
\Big\}  \otimes \psi_2 +
(\psi_1^+ + \psi_1^-)  \otimes (\nabla^2_{\pi_2 (X)} \, \psi_2 )  ,        \nonumber   \\
&    &       \\
D ( (\psi_1^+ + \psi_1^- ) \otimes \psi_2) & = &  ( D_1 \psi_1^+ +
D_1 \psi_1^-   )  \otimes \psi_2  +
 \Big\{  E_{2m+2}  \cdot  ( \psi_1^-  -  \psi_1^+ )  \Big\}  \otimes (D_2 \psi_2) ,    \nonumber    \\
&    &       \\
D^2 ((\psi_1^+ + \psi_1^-) \otimes \psi_2) & = & ( D_1^2 \psi_1^+
+ D_1^2 \psi_1^- )   \otimes \psi_2 +  ( \psi_1^+ + \psi_1^- )
\otimes ( D_2^2 \psi_2 )    ,
\end{eqnarray}
where \, $\pi_1 : T(M_1^{2m+1} \times M_2^{2s+1}) \longrightarrow
T(M_1^{2m+1})$ and $\pi_2 : T(M_1^{2m+1} \times M_2^{2s+1})
\longrightarrow T(M_2^{2s+1})$ are the natural projections.

\noindent
\begin{thm}
Let $(M_1^{2m+1}, g_1), \,  m \geq 1,$ be a Riemannian spin
manifold admitting a Killing spinor $\psi_1^{\pm} \in \Gamma (
\Sigma (M_1^{2m+1})) \subset \Gamma ( \Sigma^{\pm} ( Q^{2m+2} ))$
with $D_1 \psi_1^{\pm} = \lambda_1 \psi_1^{\pm},  \,  \lambda_1
\not=0  \in  {\mathbb R}.$   Let $(M_2^{2m+1}, g_2)$ be a
Riemannian spin manifold admitting a Killing spinor $\psi_2$ with
$D_2 \psi_2 = \lambda_2, \, \lambda_2  \in {\mathbb R}$ (Here we
allow $\lambda_2$ to be zero).  Then the Riemannian product
manifold $(M_1^{2m+1} \times M_2^{2m+1} , \, g_1 + g_2)$ admits
such a non-trivial eigenspinor $\psi^{\ast}$ of the Dirac operator
$D$ that satisfies the equation (1.3) of Theorem 1.1.
\end{thm}

\noindent {\bf Proof.}
 Set
\[    \psi : = ( \psi_1^+ + \psi_1^-) \otimes \psi_2 .     \]
Then, using the formula (4.11), we have
\[    D^2 \psi  =  ( \lambda^2_1 + \lambda^2_2 ) (\psi^+_1 + \psi_1^- )  \otimes \psi_2  =  ( \lambda^2_1 + \lambda^2_2 )   \psi   \]
and  $\lambda^{\ast} : =  \sqrt{\lambda^2_1 + \lambda^2_2}$  is
related to the scalar curvature $S = S_1 + S_2$ of $(M_1^{2m+1}
\times M_2^{2m+1} , \, g_1 + g_2 )$ as
\[   (\lambda^{\ast})^2 = \frac{1}{4} \cdot   \frac{2m+1}{2m+1 -1}  \cdot  S_1   +  \frac{1}{4} \cdot   \frac{2m+1}{2m+1-1}  \cdot  S_2   =
\frac{1}{4}  \cdot  \frac{n}{n-2}  \cdot  S ,  \quad  n = 4m+2.
\] Consequently,
\[
\psi^{\ast}    :=    \lambda^{\ast} \, \psi + D \psi  =  (
\lambda^{\ast}  + \lambda_1  )  ( \psi_1^+ + \psi_1^- ) \otimes
\psi_2
  +  \lambda_2   \{  E_{2m+2}  \cdot  ( \psi^- -  \psi^+ )  \}  \otimes  \psi_2
 \]
is indeed such a non-trivial eigenspinor of $D$  with eigenvalue
$\lambda^{\ast}$ that satisfies the equation (1.3). \hfill{QED.}

\bigskip  \noindent
{\bf  Remark.}   As mentioned above,   the spinor  bundle
$\Sigma(M_1^{2m+1})$ of $(M_1^{2m+1} ,   g_1)$ may be thought to
be embedded in the positive part $\Sigma^+( Q^{2m+2} )$ (resp. in
the negative part $\Sigma^- (Q^{2m+2})$ ) of the spinor bundle
$\Sigma( Q^{2m+2} )$.  Obviously,  there exists  a Killing spinor
$\psi_1^+  \in \Gamma ( \Sigma (M_1^{2m+1})) \subset \Gamma (
\Sigma^+ ( Q^{2m+2} ))$ with $D_1 \psi_1^+ = \lambda_1 \psi_1^+$
if and only if
 there exists  a Killing spinor $\psi_1^-  \in \Gamma ( \Sigma (M_1^{2m+1})) \subset \Gamma ( \Sigma^- ( Q^{2m+2} ))$ with $D_1 \psi_1^- = \lambda_1 \psi_1^-$ ( e.g., one can take $\psi_1^- :=  E_{2m+2}   \cdot  \psi_1^+$  ).

\noindent
\begin{thm}
Let $(M_1^{2m+1}, g_1), \,  m \geq 1,$ be a Riemannian spin
manifold admitting a Killing spinor $\psi_1^{\pm} \in \Gamma (
\Sigma (M_1^{2m+1})) \subset \Gamma ( \Sigma^{\pm} ( Q^{2m+2} ))$
with $D_1 \psi_1^{\pm} = \lambda_1 \psi_1^{\pm}, \, \lambda_1
\not=0  \in  {\mathbb R} .$ Let $( \overline{M}_2^{\, 2s+1},
\overline{g}_2), \, 2m+1  >  2s +1  \geq  3 ,$ be a Riemannian
spin manifold admitting a parallel spinor $\psi_2$. Let $( S^1,
g_S)$ be a circle with the standard metric. Then the Riemannian
product manifold $(M_1^{2m+1} \times  \overline{M}_2^{\, 2s+1} ,
\, g_1 +  \overline{g}_2)$ as well as
 $( M_1^{2m + 1}  \times S^1, \,
g_1 + g_S )$  admits  such a non-trivial eigenspinor $\psi^{\ast}$
of the Dirac operator $D$ that satisfies the equation (1.4) of
Theorem 1.2
\end{thm}

\noindent {\bf  Proof.}   Applying the argument  in the proof  for
Theorem 4.3,  one proves that $(M_1^{2m+1} \times
\overline{M}_2^{\, 2s+1} , \, g_1 +  \overline{g}_2)$ satisfies
the limiting case of Theorem 1.2.    To prove the latter part of
the theorem,  it suffices to check that  $\psi_1^{\ast} : =
\psi_1^+ +    \psi_1^-$   satisfies the equation (1.4) of Theorem
1.2.   \hfill{QED.}

\bigskip  \noindent
{\bf Remark.}  Let $(M^n, g), \, n \geq 3,$ be a Riemannian spin
manifold possessing a parallel unit vector field $\xi$. Let $\eta
= g( \cdot, \, \xi)$ be the dual 1-form. Then the endomorphism
$\phi$ defined by
\[   \phi(X) = X - 2 \, \eta(X) \xi  \]
is an almost product Riemannian structure with  $\nabla  \phi  =
0$.   Thus we find that the inequality (1.2) is indeed a special
case of the inequality in Theorem 1.2, with  $n_1 = n-1$ and $n_2
=1$.

\bigskip  \noindent
{\bf Remark.} It may be of interest to classify  all the types of
limiting manifolds of Theorem 1.1 (resp. Theorem 1.2). An
important problem  toward this classification is to consider only
simply connected  limiting manifolds, i.e., those limiting
manifolds $(M^n, g)$ that are global Riemannian products $(M^n, g)
= (M_1^{n_1}  \times M_2^{n_2},  \, g_1 + g_2 )$, and  answer the
following question : Do there exist such limiting manifolds
(Riemannian products) of Theorem 1.1 (resp. Theorem 1.2) that do
not belong to the type (i) (resp. the type (ii))  described at the
end of the introduction of the paper ?

\bigskip \bigskip  \noindent
{\bf Acknowledgements.}   This research was supported by the
Special Research Fund 2003 of  Andong National University.

\begin{Literature}{xx}
\bibitem{1}
Alexandrov, B., Grantcharov, G., Ivanov, S.:  An estimate for the
first eigenvalue of the Dirac operator on compact Riemannian spin
manifold admitting parallel one-form. J. Geom. Phys. {\bf 28},
263-270 (1998).
\bibitem{2}
B\"{a}r, C.:  Real Killing spinors and holonomy.  Commun. Math.
Phys. {\bf 154}, 509-521  (1993)
\bibitem{3}
Baum, H.,  Friedrich, Th., Grunewald, R.,  Kath, I.:  Twistors and
Killing spinors on Riemannian manifolds. Leipzig-Stuttgart:
Teubner 1991
\bibitem{4}
Cahen, M., Gutt, S., Trautman  A.:  Pin structures and the
modified  Dirac operator.  J.  Geom. Phys.  {\bf 17},  283-297
(1995)
\bibitem{5}
Hijazi, O.:   A conformal lower bound for the smallest eigenvalue
of the Dirac operator and Killing spinors.  Commun. Math. Phys.
{\bf 104},  151-162 (1986)
\bibitem{6}
Hijazi, O.:   Lower bounds for eigenvalues of the Dirac operator
through modified connections.   J. Geom. Phys. {\bf 16},  27-38
(1995)
\bibitem{7}
Hijazi, O.:  Spectral properties of the Dirac operator and
geometrical structures. Univ. Nantes: Lecture note 1998
\bibitem{8}
Friedrich,  Th.:  Der erste Eigenwert des Dirac-Operators einer
kompakten Riemannschen Mannigfaltigkeit nichtnegativer
Skalarkr\"{u}mmung.  Math. Nachr. {\bf 97},  117-146 (1980)
\bibitem{9}
Friedrich, Th.: The classification of 4-dimensional K\"{a}hler
manifolds with small eigenvalue of the Dirac operator. Math. Ann.
{\bf 295},  565-574 (1993)
\bibitem{10}
Friedrich, Th., Kim, E.C.:  The Einstein-Dirac equation on
Riemannian spin manifolds.  J. Geom. Phys. {\bf 33}, 128-172
(2000)
\bibitem{11}
Friedrich, Th., Kim, E.C.:   Some remarks on the Hijazi inequality
and generalizations of the Killing equation for spinors.  J. Geom.
Phys. {\bf 37},  1-14 (2001)
\bibitem{12}
Jung, S.D.: The first eigenvalue of the transversal Dirac
operator.  J. Geom. Phys. {\bf 39},  253-264 (2001)

\bibitem{13}
Kirchberg, K.-D.:  An estimation for the first eigenvalue of the
Dirac operator on closed K\"{a}hler manifolds of positive scalar
curvature. Ann. Global Anal. Geom. {\bf 4},  291-326 (1986)
\bibitem{14}
Kirchberg, K.-D.:   The first eigenvalue of the Dirac operator on
K\"{a}hler manifolds.  J. Geom. Phys. {\bf 7},  449-468 (1990)
\bibitem{15}
Kirchberg, K.-D.:    Properties of K\"{a}hlerian twistor-spinors
and vanishing theorems.  Math. Ann. {\bf 293}, 349-369 (1992)
\bibitem{16}
Kramer, W., Semmelmann, U., Weingart, G.:   Eigenvalue estimates
for the Dirac operator on quaternionic K\"{a}hler manifolds. Math.
Z. {\bf 230},  727-751 (1999)
\bibitem{17}
Kraus, M., Tretter, C.:  A new method for eigenvalue estimates for
Dirac operators on certain manifolds with $S^k$-symmetry.  Differ.
Geom. Appl. {\bf 19},  1-14 (2003)
\bibitem{18}
Lott, J.:   Eigenvalue bounds for the Dirac operator.  Pac. J.
Math. {\bf 125},  117-126 (1986)
\bibitem{19}
Moroianu, A.:   K\"{a}hler manifolds with small eigenvalues of the
Dirac operator and a conjecture of Lichnerowicz. Ann. Inst.
Fourier (Grenoble)  {\bf 49}, 1637-1659  (1999)
\bibitem{20}
Moroianu, A., Ornea, L.:   Eigenvalue estimates for the Dirac
operator and harmonic 1-forms of constant length. math.DG/0305140.
\bibitem{21}
Petersen, P.:    Riemannian geometry. New York-Berlin-Heidelberg:
Springer 1998
\bibitem{22}
Wang, M.:   Parallel spinors and parallel forms.  Ann. Global
Anal. Geom. {\bf 7}, 59-68  (1989)
\bibitem{23}
Yano, K., Kon, M.:  Structure on manifolds.  Singapore: World Sci.
1984
\end{Literature}

\end{document}